\documentclass{AIMS}
\usepackage{amsmath, amssymb}
\usepackage{lmodern}
\usepackage{microtype}

\usepackage{paralist}
\usepackage{graphics} 
\usepackage{epsfig} 
\usepackage{graphicx}  \usepackage{epstopdf}
 \usepackage[colorlinks=true]{hyperref}
\hypersetup{urlcolor=blue, citecolor=red}

\usepackage[numbers,sort&compress]{natbib}
\usepackage{bm}

\newcommand{\mat}[1]{\bm{#1}}
\renewcommand{\vec}[1]{\bm{#1}}
\newcommand{\koopman}{\mathcal{K}}
\newcommand{\manifold}{\mathcal{M}}
\newcommand{\evolution}{\vec{F}}
\newcommand{\featurespace}{\mathcal{F}}
\newcommand{\mode}{\vec\xi}
\newcommand{\range}{\mathcal{R}}

\def\currentvolume{X}
 \def\currentissue{X}
  \def\currentyear{200X}
   \def\currentmonth{XX}
    \def\ppages{X--XX}
     \def\DOI{10.3934/xx.xx.xx.xx}

\theoremstyle{definition}

\title[Kernel-Based Methods for Koopman Spectral Analysis] 
      {A Kernel-Based Method \\ for Data-Driven Koopman Spectral Analysis}

\author[M.O. Williams \and C.W. Rowley \and I.G. Kevrekidis]{}


 \email{mow2@princeton.edu}
 \email{cwrowley@princeton.edu}
 \email{yannis@princeton.edu}

\begin{document}
\maketitle

\centerline{\scshape Matthew O. Williams }
\medskip
{\footnotesize
 \centerline{Program in Applied and Computational Mathematics (PACM)}
 \centerline{Princeton University}
 \centerline{Princeton, NJ 08544, USA}
}

\medskip

\centerline{\scshape Clarence W. Rowley}
\medskip
{\footnotesize
\centerline{Department of Mechanical and Aerospace Engineering}
\centerline{Princeton University}
\centerline{Princeton, NJ 08544, USA}
}

\medskip

\centerline{\scshape Ioannis G. Kevrekidis}
\medskip
{\footnotesize
 \centerline{Department of Chemical and Biological Engineering \& PACM}
 \centerline{Princeton University}
 \centerline{Princeton, NJ 08544, USA}
}

\bigskip

 \centerline{(Communicated by the associate editor name)}

\begin{abstract}
A data-driven, kernel-based method for approximating the leading Koopman
eigenvalues, eigenfunctions, and modes in problems with high-dimensional state
spaces is presented.
This approach approximates the Koopman operator using a set of scalar
observables, which are functions that map states to scalars, that is determined {\em implicitly\/} by the choice of a kernel function.
This circumvents the computational issues that arise due to the number of basis
functions required to span a ``sufficiently rich'' subspace of the space
of scalar observables in such applications.
We illustrate this method on two examples: the FitzHugh-Nagumo PDE,
 a prototypical one-dimensional reaction-diffusion
system, and vorticity data obtained from experimentally obtained
velocity data for flow over a cylinder at Reynolds number 413.
In both examples, we compare our  results
with related methods, such as Dynamic Mode Decomposition (DMD) that
have the same cost as our approach.
\end{abstract}

\section{Introduction}

In many applications, the evolution of complex spatio-temporal phenomena  can
be characterized using models based on the interactions between a relatively
small number of modes.
Due to the availability of data and computational power, algorithmic techniques
for identifying such modes have become increasingly
common~\cite{Holmes1998,chatterjee2000introduction,Schmid2011,Schmid2012,Rowley2009,Budivsic2012,WilliamsEDMD,Wynn2013,Juang1994,Chen2012,Tu2013,Tissot2014,Coifman2006,Lee2007}.
Perhaps the best known method is the Proper Orthogonal Decomposition
(POD)~\cite{Holmes1998,chatterjee2000introduction,Bishop2006}, which is
also known as Principal Component Analysis (PCA)~\cite{Bishop2006,Lee2007}.
However, other algorithms that generate different sets of modes exist.
In recent years, approximations of the modes of the Koopman
operator~\cite{Koopman1931,Koopman1932} have become popular in fluid
applications~\cite{Rowley2009,Budivsic2012,mezic2013analysis,Bagheri2013,Chen2012,Wynn2013,Hemati2014}.
The Koopman modes of the full state observable, which are vectors for systems of
ODEs or spatial profiles for PDEs, are intrinsic to a particular evolution law,
and have temporal dynamics that are determined by their corresponding Koopman
eigenvalues~\cite{Rowley2009};
in effect, the Koopman modes look and act like the eigenvectors of a linear
system even when the underlying evolution law is nonlinear.

This representation is possible because the Koopman operator, which
defines these quantities, is a {\em linear\/} operator and therefore can have
eigenfunctions and eigenvalues, which are later used to define the modes.
However, it acts on {\em scalar observables}, which are functions defined on
state space, and is therefore an infinite-dimensional operator even for finite-dimensional
dynamical systems.
As a result,  methods that approximate the Koopman operator (often
implicitly) select a finite-dimensional subspace of the space of scalar
observables to use during the computation.
Currently, the most widely used method is Dynamic Mode Decomposition
(DMD)~\cite{Schmid2010,Schmid2011,Schmid2012,Rowley2009,Budivsic2012} (along
with its related modifications~\cite{Chen2012,Wynn2013,Hemati2014}),
which implicitly select \textit{linear} functions of the state as the basis functions~\cite{WilliamsEDMD}.
This restrictive choice of basis allows DMD to be applied to large systems of
ODEs or PDEs, but in many applications, this subspace is simply not ``rich''
enough to effectively approximate the Koopman operator~\cite{WilliamsEDMD,Tu2013}.

Other methods, such as Extended DMD~\cite{WilliamsEDMD}, use an expanded set of
observables, which can produce a more accurate approximation of the
set of Koopman
eigenvalues, eigenfunctions, and modes  for nonlinear systems.
For small systems of ODEs, a ``sufficiently rich'' subspace can often be
spanned using a few hundred or thousand basis functions, which is
computationally tractable even on a laptop.
However, the size of the needed basis grows rapidly as the dimension of the
state space increases, so the necessary
computations quickly become infeasible.
This explosion in the computational cost is common in machine learning
applications, and one facet of the ``curse of
dimensionality''~\cite{Bishop2006}.

In this manuscript, we introduce a data-driven, kernel-based method to
approximate the Koopman operator in systems with large state dimension.
This approach circumvents the dimensionality issues encountered by Extended DMD
by defining a kernel function that {\em implicitly\/} computes inner products in
the high-dimensional space of observables.
Because we do not form this space explicitly, the computational
cost of the method is determined by the number of snapshots and the
dimension of state space rather than the number of basis
functions used to represent the scalar observables;
therefore, the cost of this approach is comparable to that of DMD.
We apply this approach to a pair of sample problems: the one-dimensional
FitzHugh-Nagumo PDE, a prototypical
reaction-diffusion system; and
experimentally obtained vorticity data from flow over a cylinder
obtained at Reynolds number 413.
Using the FitzHugh-Nagumo PDE example, we will demonstrate that using
a higher-dimensional subspace of observables can result in more accurate and
reproducible approximations of the leading Koopman modes, eigenvalues, and eigenfunctions.
Finally, the cylinder example demonstrates that this approach has practical
advantages in more realistic settings where the data are noisy and the
``true'' Koopman modes and eigenvalues are unknown.

The remainder of the manuscript is outlined as follows: in
Sec.~\ref{sec:numerical-method} we present the kernel reformulation of
Extended Dynamic Mode Decomposition.
In Sec.~\ref{sec:exampl-fitzh-nagumo}, we apply it to a numerical
discretization of the FitzHugh-Nagumo PDE in one spatial dimension,
where a subset of the Koopman eigenvalues and modes are known.
In Sec.~\ref{sec:exampl-exper-data}, we apply our method to
experimentally obtained data of flow over a cylinder, which is a more
realistic example with measurement noise and other experimental
realities.
Finally, in Sec.~\ref{sec:conclusions}, we present some concluding
remarks and future outlook. 

\section{A Data-Driven Approximation of the Koopman Operator}
\label{sec:numerical-method}

In this section, we present a data-driven approach for approximating
the Koopman operator that can be applied to systems with
high-dimensional state spaces. 
The method presented here is a reformulation of the Extended DMD
procedure that makes use of the so-called {\em kernel trick}~\cite{Scholkopf2001}.
In the following subsections, we will (i) briefly review the Koopman
operator (ii) give a brief derivation
of ``standard'' Extended DMD, which enables the ``tuples'' of Koopman
eigenvalues, eigenfunctions, and modes to be approximated from data; (iii)
present the kernel approach; and  (iv) give a practical algorithm for computing the leading tuples of eigenvalues, eigenfunctions, and modes.

\subsection{The Koopman Operator}
\label{sec:koopman}

The Koopman operator~\cite{Koopman1931,Koopman1932}, along with its eigenvalues, eigenfunctions, and
modes, is defined by a dynamical system and not a set of data.
Given the discrete time dynamical system $(n,\mathcal{M},\vec F)$,
where $n\in\mathbb{Z}$ is  time, $\mathcal{M}\subseteq \mathbb{R}^N$ is the state space, and $\vec x
\mapsto \vec F(\vec x)$ defines the dynamics, the action of the Koopman operator, $\mathcal{K}$, on a
{\em scalar observable}, $\phi:\mathcal{M}\to\mathbb{C}$, is
\begin{equation}
(\mathcal{K}\phi)(\vec x) = (\phi\circ \vec F)(\vec x) = \phi(\vec F(\vec x)).
\end{equation}
Intuitively,  $\mathcal{K}\phi$ is a new function that gives the value of
$\phi$ ``one step in the future''.
Note that {\em the Koopman operator acts on functions of the state, and not the states themselves}.
Since $\phi\in\mathcal{F}$, where $\mathcal{F}$ is an appropriate space of
scalar observables, the Koopman operator is infinite dimensional.
However, it is also linear, and thus it can
have eigenvalues and eigenfunctions, which we refer to as Koopman eigenvalues
and eigenfunctions.
Accompanying the eigenvalues and eigenfunctions are the Koopman modes for a
given {\em vector valued observable}, $\vec g:\mathcal{M}\to\mathbb{R}^{N_o}$, where $N_o\in\mathbb{N}$.
These modes are vectors in a system of ODEs (or spatial profiles in a PDE) that
contain the coefficients required to construct $\vec g$ using a Koopman
eigenfunction basis~\cite{Rowley2009,WilliamsEDMD}.
One particularly useful set of modes is that of the identity operator or, equivalently, the {\em full state
observable},  $\vec g(\vec x) = \vec x$,  which we refer to as simply the
{\em Koopman modes\/} in all that follows.

In many systems~\cite{Rowley2009,Mezic2005,mezic2013analysis},  tuples
consisting of an eigenvalue, an eigenfunction, and a mode enable a simple yet
powerful means of representing the system state and making predictions of future
values.
In particular,
\begin{equation}
	\vec x = \sum_{k=1}^K \vec\xi_k \varphi_k(\vec x),\qquad
	\vec F(\vec x) = \sum_{k=1}^K \mu_k\vec\xi_k \varphi_k(\vec x)
	\label{eq:koopman_modes}
\end{equation}
where $\vec\xi_k$ is the Koopman mode corresponding to the eigenfunction
$\varphi_k$, $\mu_k$ is the corresponding eigenvalue, and $K$ is the number of tuples
required for the reconstruction, which could be infinite.
The eigenfunction $\varphi_k:\mathcal{M}\to\mathbb{C}$ is a {\em
function}, but the mode $\vec\xi_k\in\mathbb{C}^N$ is a {\em vector}.
From~\eqref{eq:koopman_modes}, the coefficient associated with the $k$-th mode,
$\vec\xi_k$,  is obtained by evaluating the $k$-th eigenfunction, $\varphi_k$,
and the temporal evolution is dictated by $\mu_k$.
Ultimately, one benefit of this Koopman-based approach is that  the
eigenvalues, eigenfunctions, and modes are {\em intrinsic to the dynamical
system}; in the subsequent section they will be approximated using data, but unlike the Proper Orthogonal Decomposition (POD) and
related methods, they still exist absent any data.

\subsection{Extended Dynamic Mode Decomposition}
\label{sec:extended-dmd}

Extended Dynamic Mode Decomposition (Extended DMD)~\cite{WilliamsEDMD}
is a regression procedure whose solution produces a
finite-dimensional approximation of the Koopman operator.
To obtain this approximation, we define
a basis set that consists of $K$ scalar observables, which we denote as
$\psi_k$ for $k=1,\ldots,K$,  that span
$\featurespace_K\subset\featurespace$.
We also define the vector valued observable,
$\vec\psi:\manifold\to\mathbb{R}^K$, where 
\begin{equation}
\vec\psi(\vec x) = 
	\begin{bmatrix}
	\psi_1(\vec x) \\ 
	\psi_2(\vec x) \\
	\vdots \\
	\psi_K(\vec x)
	\end{bmatrix}.
\end{equation}
In this application, $\vec\psi$ is the mapping from
physical space to {\em feature space}.
Any $\phi,\hat\phi \in\featurespace_K$ can be written as 
\begin{equation}
\phi = \sum_{k=1}^K a_k\psi_k = \vec\psi^T\vec a, \qquad \hat\phi = \sum_{k=1}^K
\hat a_k\psi_k =\vec\psi^T\hat{\vec a},
\label{eq:finite}
\end{equation}
for some set of coefficients $\vec a,\hat{\vec a}\in\mathbb{C}^K$.
Although $\evolution$ is unknown, we assume access to a data set of
snapshot \textit{pairs}:
\begin{equation}
\{(\vec x_m, \vec y_m)\}_{m=1}^M,
\qquad \vec y_m = \evolution(\vec x_m),
\end{equation}
where $\vec x_m,\vec y_m\in\manifold$.
One important special case of such a data set is a single time-series of data,
which can be written in the above form by ``grouping'' sequential
pairs of snapshots that were obtained with a fixed sampling interval,
say $\Delta t$.

Let $\hat\phi = \koopman\phi + r$, where $r\in\featurespace$ is a residual
function that appears because $\featurespace_K$ is not necessarily invariant to
the action of the Koopman operator.
Using the notation in \eqref{eq:finite}, the objective of the Extended DMD
procedure~\cite{WilliamsEDMD} is to define a mapping from some given $\vec a$ to
a  new vector $\hat{\vec a}$ that minimizes this residual.
Because the Koopman operator is linear, this mapping can be represented by a
matrix $\mat{K}\in\mathbb{R}^{K\times K}$.
To determine the entries of $\mat{K}$, the Extended DMD approach takes ideas
from collocation methods typically used to solve PDEs, but uses the $\vec x_m$
as the collocation points rather than a pre-determined
grid~\cite{Boyd2013,Trefethen2000}.
As a result, the finite dimensional approximation is 
\begin{subequations}
\label{eq:edmd-approximation}
\begin{equation}
\mat{K} \triangleq \mat{\Psi}_x^+\mat{\Psi}_y,
\end{equation}
where 
\begin{equation}
\mat{\Psi}_x \triangleq 
	\begin{bmatrix}
		\vec\psi(\vec x_1)^T \\
		\vec\psi(\vec x_2)^T \\
		\vdots \\
		\vec\psi(\vec x_M)^T
	\end{bmatrix}, \qquad 
\mat{\Psi}_y \triangleq 
	\begin{bmatrix}	
		\vec\psi(\vec y_1)^T \\
		\vec\psi(\vec y_2)^T \\
		\vdots \\
		\vec\psi(\vec y_M)^T
	\end{bmatrix},
\end{equation}
\end{subequations}
are in $\mathbb{R}^{M\times K}$, and $+$ denotes the pseudoinverse.

In \citet{WilliamsEDMD}, the relationship $\mat{\Psi}_x^+ = (\mat{\Psi}_x^T\mat{\Psi}_x)^+\mat{\Psi}_x^T$ was used to rewrite \eqref{eq:edmd-approximation} as 
\begin{subequations}
  \label{eq:edmd-normal}
  \begin{equation}
    \mat{K} \triangleq \mat{G}^+\mat{A},
  \end{equation}
  \begin{equation}
    \mat{G} =  \mat{\Psi}_x^T\mat{\Psi}_x = \sum_{m=1}^M \vec\psi(\vec x_m)\vec\psi(\vec x_m)^T, \qquad 
    \mat{A} =  \mat{\Psi}_x^T\mat{\Psi}_y = \sum_{m=1}^M \vec\psi(\vec x_m)\vec\psi(\vec y_m)^T,     
  \end{equation}
\end{subequations}
 which is advantageous when the number of snapshots is much larger than the number of basis functions, i.e., $M\gg K$, because each term in the sum can be evaluated individually.
As a result, one need only store the matrices $\mat{G},\mat{A}\in\mathbb{R}^{K\times K}$ rather than $\mat{\Psi}_x,\mat{\Psi}_y\in\mathbb{R}^{M\times K}$, which would be much larger in this regime. 

Expressions \eqref{eq:edmd-approximation} and~\eqref{eq:edmd-normal} are mathematically equivalent, and produce the same matrix~$\mat{K}$.
Therefore, regardless of how it was computed, the properties of $\mat{K}$ of interest here are unchanged:
\begin{enumerate}
\item The $k$-th eigenvalue of $\mat{K}$, $\mu_k$, is an approximation
  of an eigenvalue of $\koopman$.
  When the data are generated by sampling a continuous time dynamical
  system with a fixed sampling interval $\Delta t$, we also define the
  approximation of the {\em continuous time eigenvalue} as $\lambda_k
  \triangleq   \log(\mu_k)/\Delta t$.
\item Using (\ref{eq:finite}), the corresponding
  eigenvector, $\vec v_k$, approximates an {\em eigenfunction} of the
  Koopman operator via 
  \begin{equation}
    \label{eq:eigenfunction-approx}
    \varphi_k \triangleq \vec\psi^T\vec v_k.
  \end{equation}
\item The left eigenvector, $\vec w_k$, can be used to approximate
  the Koopman mode, $\vec\xi_k$.
\end{enumerate}
For a derivation of these relationships, see \citet{WilliamsEDMD}.

\subsection{The Kernel Method}
\label{sec:kernel}

In Sec.~\ref{sec:extended-dmd}, we considered the case where the number of snapshots was large compared to the dimension of the desired subspace of functions.
Now we will consider the opposite and more commonly encountered regime~\cite{Schmid2012,Schmid2011,Schmid2010,Rowley2009,Budivsic2012,mezic2013analysis,Bagheri2013,Chen2012,Wynn2013,Hemati2014}  where the number of snapshots is small compared to the dimension of our subspace of scalar observables (i.e., $M \ll  K$).

The difficulty here is that the  Extended DMD procedure, as formulated in Sec.~\ref{sec:extended-dmd}, requires a $K \times K$ matrix to be formed and decomposed, which requires $\mathcal{O}(K^2M)$ and $\mathcal{O}(K^3)$ time respectively, and the value of $K$ for a ``rich'' set of basis functions grows rapidly as the dimension of state space increases.
For instance, consider the case where  $\mathcal{F}_K$ is the space of all (multivariate) polynomials on
$\mathbb{R}^{256}$ with degree up to 20, as it will be in our
first example discussed in Section~\ref{sec:exampl-fitzh-nagumo}.
In this case, $K\sim 10^{30}$, which is far too large for practical computations~\cite{Bishop2006}.
This explosion in the size of the problem is common, and is one facet of
the {\em curse of dimensionality}.

Because the matrix $\mat{K}$ is the solution to a regression problem, the non-zero eigenvalues and their associated left and right eigenvectors can also be obtained by solving the {\em dual form}  of this problem~\cite{Bishop2006}.
To show this, note that $\range(\mat{K})\subseteq\range(\mat{\Psi}_x^T)$, i.e., the range of
$\mat{\Psi}_x^T$  contains the range of $\mat{K}$.
If we could compute the SVD of $\mat{\Psi}_x$,
\begin{equation}
\mat{\Psi}_x \triangleq \mat{Q}\mat{\Sigma}\mat{Z}^T,
\label{eq:psi-svd}
\end{equation}
where $\mat{Q},\mat{\Sigma}\in\mathbb{R}^{M\times M}$ and $\mat{Z} \in\mathbb{R}^{K\times M}$,  then an eigenvector of $\mat{K}$ with $\mu_k\neq 0$ could be written as $\vec v = \mat{Z}\vec{\hat{v}}$ for some $\vec{\hat{v}}\in\mathbb{C}^M$.
With some simple algebraic manipulations, the eigenvalue problem can be written as:
\begin{align*}
\mu\vec{v} &= \mat{K}\vec v\\
\iff \mu\mat{Z}\hat{\vec v}
  &= \mat{Z}\mat{\Sigma}^+\mat{Q}^T\mat{\Psi}_y\mat{Z}\vec{\hat{v}} \\
  &= \mat{Z}\left(\mat{\Sigma}^+\mat{Q}^T\right)\left(\mat{\Psi}_y\mat{Z}\right)\vec{\hat{v}} \\ 
  &= 
 \mat{Z}\left[\left(\mat{\Sigma}^+\mat{Q}^T\right)\left(\mat{\Psi}_y\mat{\Psi}_x^T\right)\left(\mat{Q}\mat{\Sigma}^+\right)\right]\vec{\hat{v}}.
\end{align*}
Therefore, an alternative method for computing an eigenvector of $\mat{K}$ is to form the matrix
 \begin{equation}
 \mat{\hat{K}} \triangleq  \left(\mat{\Sigma}^+\mat{Q}^T\right)\mat{\hat{A}}\left(\mat{Q}\mat{\Sigma}^+\right),
 \label{eq:kernel-edmd}
 \end{equation}
 where $\mat{\hat{A}} \triangleq \mat{\Psi}_y\mat{\Psi}_x^T$, 
compute an eigenvector of $\mat{\hat{K}}$, say $\hat{\vec v}$, and set $\vec v = \mat{Z}\vec{\hat{v}}$.
Here $\mat{\hat{K}}\in\mathbb{R}^{M\times M}$, so the computational cost
of the decomposition is determined \textit{by the number of snapshots} rather than the
dimension of the system state or ``feature'' space.

The benefit of the expression in $\mat{\hat{K}}$ is that all the required matrices can be obtained by computing inner products in feature space. 
In addition to $\mat{\hat{A}}$, we define the matrix $\mat{\hat{G}} \triangleq \mat{\Psi}_x\mat{\Psi}_x^T$.
Despite appearing ``flipped,'' the $ij$-th elements of $\mat{\hat{G}}$ and $\mat{\hat{A}}$ are
\begin{align}
  \label{eq:entries}
  \mat{\hat{G}}^{(ij)} \triangleq \vec\psi(\vec{x}_j)^T\vec\psi(\vec{x}_i),\qquad
  \mat{\hat{A}}^{(ij)} \triangleq \vec\psi(\vec{x}_j)^T\vec\psi(\vec{y}_i).
\end{align}
However $\mat{\hat{G}} = \mat{Q}\mat{\Sigma}^2\mat{Q}^T$, using the definition of $\mat{Q}$ and $\mat{\Sigma}$ in \eqref{eq:psi-svd}.
Therefore, given $\mat{\hat{G}}$ we can obtain $\mat{Q}$ and $\mat{\Sigma}$ via its eigen-decomposition, which is simply the {\em method of snapshots} approach for computing the POD modes from snapshot data~\cite{Sirovich1987}.
As a result, we {\em could} compute $\mat{\hat{K}}$ by forming $\mat{\hat{G}}$ and $\mat{\hat{A}}$ using \eqref{eq:entries} in $\mathcal{O}(M^2K)$ time.
This is a large improvement over ``standard'' Extended DMD, but still impractical as $K$ can be extremely large.

Rather than {\em explicitly defining} the function $\vec\psi$ and
computing the entries of $\mat{\hat{K}}$ directly, the {\em kernel trick} is a 
common technique for {\em implicitly} computing inner products~\cite{Burges1998,Baudat2001,Scholkopf2001,rasmussen2006gaussian},
which can be used to assemble $\mat{\hat{K}}$ in $\mathcal{O}(M^2N)$ time.
Instead of defining $\vec\psi$, we define a {\em kernel function\/}
${f:\mathcal{M}\times\mathcal{M}\to\mathbb{R}}$ that computes inner products in
feature space given pairs of data points; that is, $f(\vec x_i, \vec x_j) =
\left\langle \vec x_i, \vec x_j\right\rangle= \vec\psi(\vec x_j)^T\vec\psi(\vec x_i)$~\cite{Burges1998}.
In effect, the choice of $f$ defines $\vec\psi$, which is equivalent to
choosing the basis in Extended DMD.
It is, however, crucial to note that {\em $f$ does not compute these inner
products directly}.
The simplest example is the polynomial kernel 
\begin{equation}
\label{eq:polynomial-kernel}
f(\vec x, \vec z) = (1 +
\vec z^T\vec x)^2
\end{equation}
with $\vec x,\vec z\in\mathbb{R}^2$, which, when expanded,
is
\begin{align*}
f(\vec{x},\vec{z}) &=(1+x_{1}z_{1}+x_{2}z_{2})^{2}\\
	&=(1+2x_{1}z_{1}+2x_{2}z_{2}+2x_{1}x_{2}z_{1}z_{2}+x_{1}^{2}z_{1}^{2}+x_{2}^{2}z_{2}^{2})\\
	&= \vec\psi(\vec z)^T\vec\psi(\vec x)
\end{align*}
if  $\vec\psi(\vec x) = [1, \sqrt{2}x_1, \sqrt{2}x_2, \sqrt{2}x_1x_2, x_1^2,
x_2^2]$.
In general, a polynomial kernel of the form 
\begin{equation}
f(\vec x, \vec z) = (1+\vec z^T\vec x)^\alpha
\end{equation}
is equivalent to a basis that can represent all polynomials up to
and including terms of degree $\alpha$, and takes only $\mathcal{O}(N)$ time to evaluate.

In the example that follows, we use a polynomial kernel with $\alpha=20$.
We selected a polynomial kernel because the Koopman eigenfunctions are often analytic in a
disk about a fixed point~\cite{Mauroy2013a}, and polynomial kernels mimic an
$\alpha$ order power-series expansion.
The specific choice of $\alpha=20$ is more {\em ad hoc}.
In general, large values of $\alpha$ use a ``richer'' set of basis
functions, but also deleteriously impact the condition number of $\mat{\hat{G}}$.
Other choices of kernels, such as Gaussian kernels (i.e., $f(\vec x,
\vec y) = \exp(-\|\vec x - \vec y\|^2/\sigma^2)$), are also used in machine learning
applications~\cite{Bishop2006, Cristianini2000,Scholkopf2001}, and may result in
better performance (both in terms of numerical conditioning and in approximating
the Koopman operator)~\cite{Ham2004}.
The optimal choice of kernel, which is equivalent to the ideal choice of
the basis set, remains an open question.

Ultimately, the procedure for obtaining an approximation of the Koopman operator is as follows:
\begin{enumerate}
\item Using the data set of snapshot pairs and the kernel function, $f$, compute the elements of $\mat{\hat{G}}$ and $\mat{\hat{A}}$ using:
\begin{equation}
\label{eq:entries-kernel}
{\mat{\hat{G}}}^{(ij)} \triangleq f(\vec x_i, \vec x_j), \quad
{\mat{\hat{A}}}^{(ij)} \triangleq f(\vec y_i, \vec x_j).
\end{equation}
\item Compute the eigendecomposition of the Gramian $\mat{\hat{G}}$ to obtain $\mat{Q}$ and $\mat{\Sigma}$.
\item Construct $\mat{\hat{K}}$ using \eqref{eq:kernel-edmd}.
\end{enumerate}
The difference between (\ref{eq:entries}) and (\ref{eq:entries-kernel}) is that the former defines $\vec{\psi}$ explicitly, and as a result, computes the set of the needed inner products in $\mathcal{O}(M^2K)$ time, while the latter defines $\vec{\psi}$ implicitly, which allows the inner products to be computed in $\mathcal{O}(M^2N)$ time.
If a linear kernel, $f(\vec x, \vec y) = \vec x^T\vec y$ is chosen, the kernel approach outlined here is identical to DMD~\cite{Tu2013}, where the subspace of observables used to approximate the Koopman eigenfunctions is chosen using the Proper Orthogonal Decomposition.
The overarching idea here is that the kernel approach chooses this subspace differently by using what is in effect {\em Kernel Principal Component Analysis}~\cite{Bishop2006}, which exploits the kernel trick to make the computation efficient when   $M \ll K$.

\subsection{Computing the Koopman Eigenvalues, Modes, and Eigenfunctions}
\label{sec:modes-eigenvalues-eigenfunctions}

In this subsection, we show how to approximate the Koopman
eigenvalues, modes, and eigenfunctions given $\mat{\hat{K}}$.
Let $\mat{\hat{V}}$ be the matrix whose columns are the eigenvectors of
$\mat{\hat{K}}$.
Then using (\ref{eq:eigenfunction-approx}), we define the matrix of eigenfunctions values:
\begin{equation}
\label{eq:eigenfunctions}
\begin{aligned}
\mat{\Phi}_x &\triangleq  \mat{\Psi}_x\mat{Z}\mat{\hat{V}} = \left(\mat{\Psi}_x\mat{\Psi}_x^T\right) \left(\mat{Q}\mat{\Sigma}^+\right)\mat{\hat{V}}  = \mat{\hat{G}}\left(\mat{Q}\mat{\Sigma}^+\right)\mat{\hat{V}} = \left(\mat{Q}\mat{\Sigma}^2\mat{Q}^T\right)\left(\mat{Q}\mat{\Sigma}^+\right)\mat{\hat{V}} \\
	&=\mat{Q}\mat{\Sigma}_r\mat{\hat{V}}
\end{aligned}
\end{equation}
where $\mat{\Sigma}_r$ is the diagonal matrix of singular values with any entry neglected in the pseudo-inverse set to zero.
The $i$-th row of $\mat{\Phi}_x$ contains the numerically computed eigenfunctions evaluated at
$\vec x_i$.
The $k$-th numerically approximated Koopman eigenfunction can also be evaluated at a new data point via:
\begin{equation}
\begin{aligned}
\varphi_k(\vec x) &= \left(\vec\psi(\vec x)^T\mat{\Psi}_x^T\right)\left(\mat{Q}\mat{\Sigma}^+\vec{\hat{v}}_k\right) \\
&= 
\begin{bmatrix}
f(\vec x, \vec x_1) & f(\vec x, \vec x_2) & \cdots & f(\vec x, \vec x_M)
\end{bmatrix}\left(\mat{Q}\mat{\Sigma}^+\vec{\hat{v}}_k\right),
\end{aligned}
\end{equation}
using the same arguments as in \eqref{eq:eigenfunctions} though without the simplifications and cancellations that occur in that case.

To compute the Koopman modes, we use (\ref{eq:koopman_modes}), which when evaluated
at each of the data points, results in the matrix equation
\begin{equation}
\mat{X} = \mat{\Phi}_x \mat{\Xi},
\end{equation}
where 
\begin{equation}
\mat{X} \triangleq 
\begin{bmatrix}
	\vec x_1^T \\
	\vec x_2^T \\
	\vdots \\
	\vec x_M^T
\end{bmatrix},\quad 
\mat{\Xi} \triangleq 
\begin{bmatrix}
	\mode_1^T \\
	\mode_2^T \\
	\vdots \\
	\mode_M^T
\end{bmatrix} = \mat{\Phi}_x^+\mat{X} =\mat{\hat{V}}^{-1}\mat{\Sigma}^+\mat{Q}^T\mat{X}.
\end{equation}
provided that $\mat{\hat{V}}$ is full rank.
In this case, 
\begin{equation}
\mat{\hat{V}}^{-1} = 
\begin{bmatrix}
\hat{\vec w}_1^* \\
\hat{\vec w}_2^* \\
\vdots \\
\hat{\vec w}_M^* 
\end{bmatrix},
\end{equation}
where $\hat{\vec w}^*$ is a left eigenvector of $\mat{\hat{K}}$ scaled so
that $\hat{\vec w}_i^*\hat{\vec v}_j = \delta_{ij}$.
This implies that the $k$-th Koopman mode, $\vec\xi_k$, is  
\begin{equation}
\vec \xi_k = (\hat{\vec w}_k^*\mat{\Sigma}^+\mat{Q}^T\mat{X})^T,
\label{eq:mode-def}
\end{equation}
and therefore approximate Koopman tuples can be obtained via the left and right
eigenvectors of $\hat{\mat K}$, and do not require $\mat{\hat{K}}$ to be
computed in its entirety. 
For the problems considered here, the complete decomposition of $\mat{\hat{K}}$
is computed, but for problems with larger numbers of snapshots, Krylov methods
could be used to compute a leading subset of its eigenvalues and
vectors~\cite{Lehoucq1998}.

\section{Example: The FitzHugh-Nagumo PDE}
\label{sec:exampl-fitzh-nagumo}

In order to highlight the effectiveness of the kernel method, we first
apply it and, as a benchmark, DMD to the FitzHugh-Nagumo
PDE~\cite{Elmer2005} in one spatial dimension.
This example is particularly useful, because a subset of the
true Koopman eigenvalues and modes can be deduced from the system
linearization.
The governing equations are:
 \begin{subequations}
 	\label{eq:fhne}
 	\begin{align}
 	\partial_{t}v & =\partial_{xx}v+v-w-v^{3},\\
 	\partial_{t}w & =\delta\partial_{xx}w+\epsilon(v-c_{1}w-c_{0}),
 	\end{align}
 \end{subequations}
where $v$ is the activator field, $w$ is the inhibitor field, $c_{0}=-0.03$, $c_{1}=2.0$, $\delta=4.0$, $\epsilon=0.02$, for $x\in[0,20]$ with Neumann boundary conditions.
Both $v$ and $w$ are approximated using a discrete cosine transform-based spectral method with 128 basis functions each.

\begin{figure}
  \centering
  \includegraphics[width=0.8\textwidth]{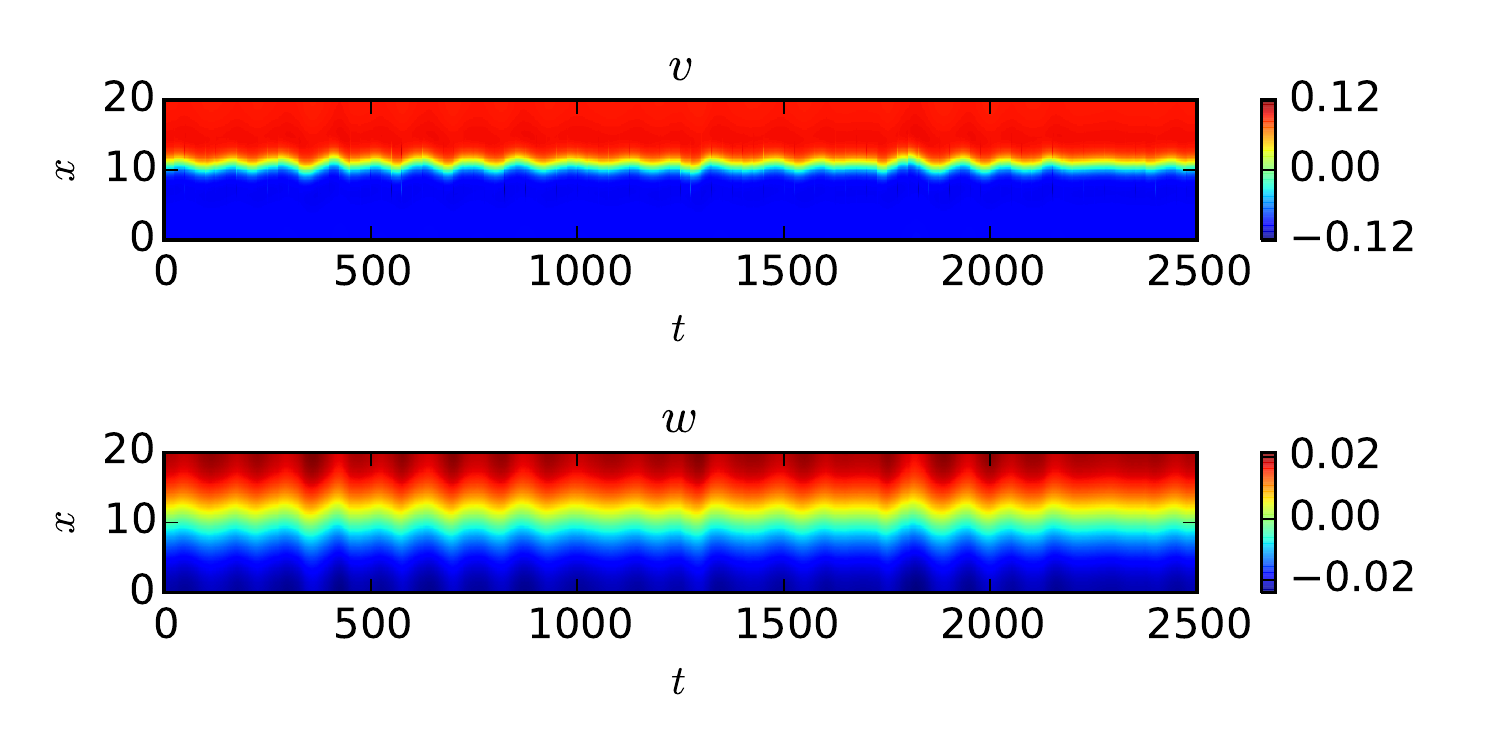}
  \caption{
  	One of the five sets of data that were used to approximate the leading Koopman modes, eigenfunctions, and eigenvalues.
  	(left) The numerically computed approximation of $v$ sampled uniformly in time and space with $\Delta t = 1$.
	(right) The $w$ data obtained at the same points in time.
	Although the initial condition for all five data sets is the equilibrium point, the random perturbations applied every 25 units of time result in ten different trajectories. 
  }
  \label{fig:fhne-data}
\end{figure}

With these parameter values, \eqref{eq:fhne} has a stable equilibrium point that resembles a standing wave-front in both the activator and inhibitor fields.
In order to explore a subset of state space and to mimic the presence of impulsive actuation, every 25 time units we perturb the system state by setting 
\begin{equation}
v(x, t) \leftarrow v(x, t) + \sum_{i=1}^3 u_i \exp(-(x-x_i)^2),
\end{equation}
where $x_1 = 7.5$, $x_2 = 10$, $x_3 = 12.5$ and  $u_1,u_2,u_3\sim \mathcal{N}(0, 0.1^2)$.
Note that $w(x,t)$ remains unchanged.
Starting from the equilibrium point, we generate five separate trajectories consisting of 2500 snapshots (or 2499 snapshot pairs) with a sampling interval of $\Delta t = 1$ starting with five different random seeds.
Figure~\ref{fig:fhne-data} shows the $v$ and $w$ data for a prototypical trajectory: the perturbations applied even 25 units in time excite ``fast'' directions that decay almost immediately (and thus are not visible in the figure) as well as ``slow'' directions that take the form of ``shifts'' in the position of the wave front.

\begin{figure}
  \centering
  \includegraphics[width=0.8\textwidth]{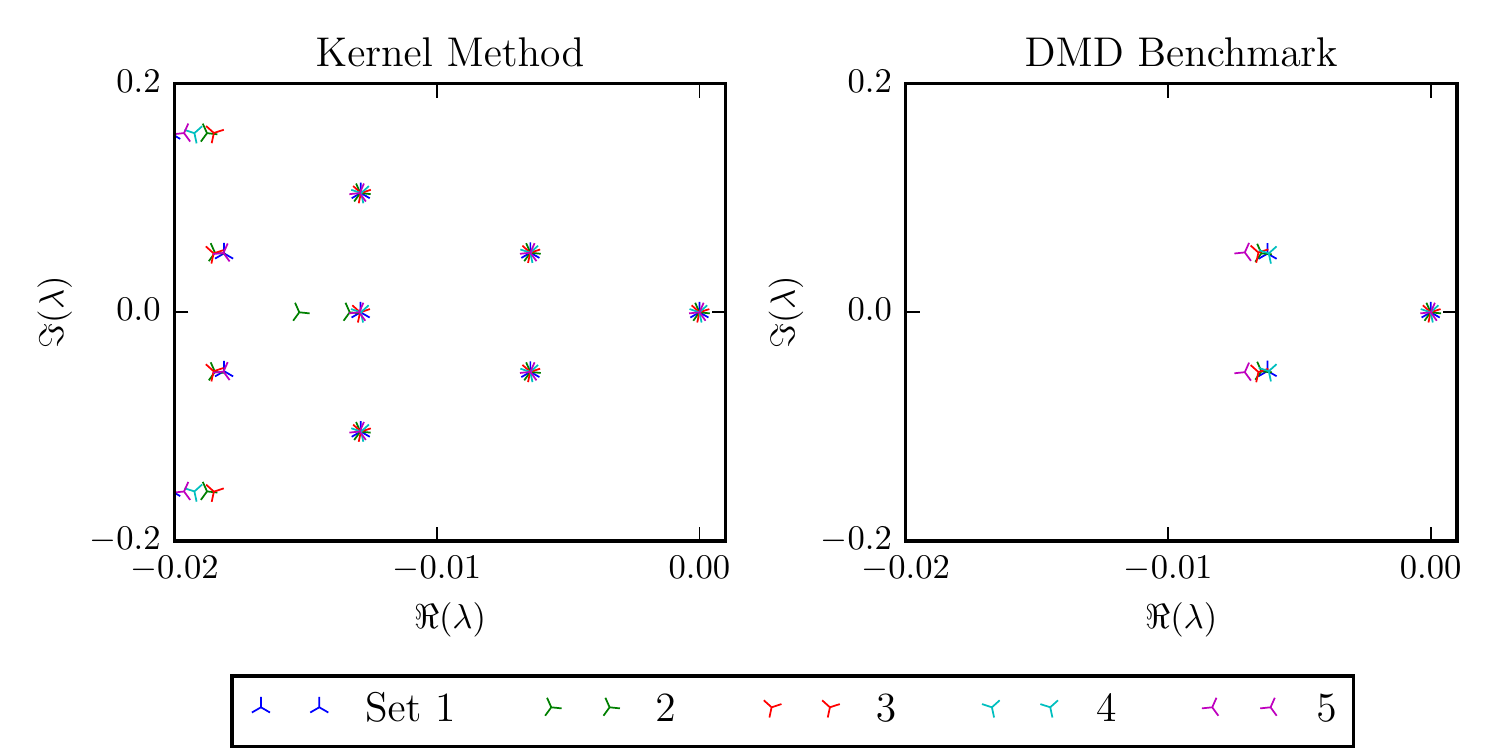}
  \caption{
    The leading eigenvalues computed using the kernel method and Dynamic Mode Decomposition are shown on the left and right respectively.
    These computations were performed on data sets like the one shown in Fig.~\ref{fig:fhne-data} using only the 150 largest singular values to compute pseudo-inverses.
    Both DMD and the kernel approach consistently identify the $\lambda_1 = 0$ eigenvalue and the pair associated with the system linearization ($\lambda_{2,3} \approx -0.006 \pm 0.053i$).
   However, the kernel-based method appears to do so with less variance, and is also able to reconstruct two additional ``layers'' of  of eigenvalues, such as $\lambda_4 \approx -0.013$ and $\lambda_7 \approx -0.019 + 0.053i$, that are known to exist analytically~\cite{Gaspard1995}.
     }
  \label{fig:fhne-eigenvalues}
\end{figure}

\begin{figure}
  \centering
  \includegraphics[width=0.8\textwidth]{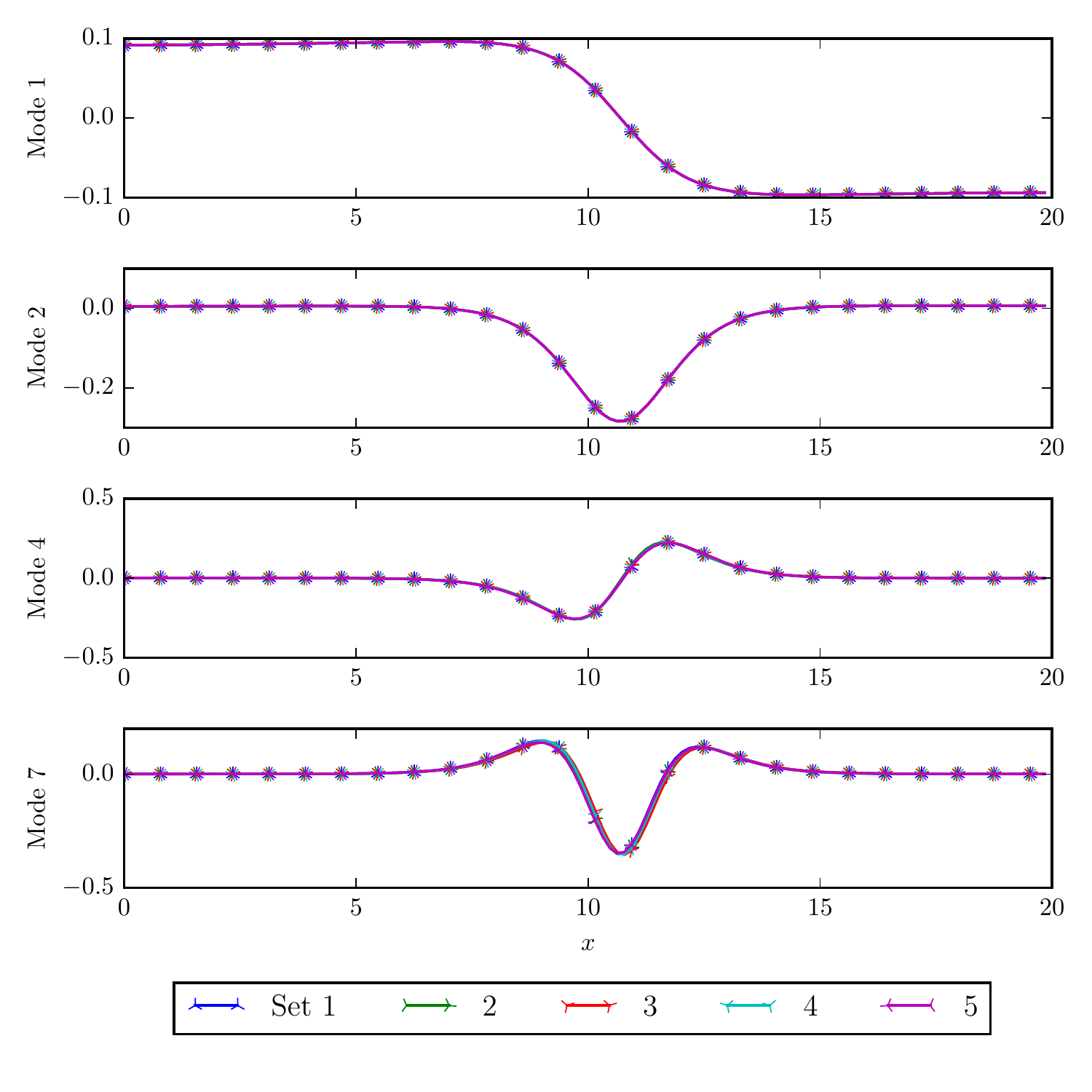}
  \caption{
    The real part of four of the Koopman modes for $v$ computed using the kernel method associated with the eigenvalues near $\lambda_1 = 0$, $\lambda_2 =  -0.006 + 0.053i$,  $\lambda_4 = -0.013$, and $\lambda_7 = -0.019 + 0.053i$ respectively. 
    Modes 1 and 2 are scaled versions of the equilibrium point and an eigenvector of the linearized system respectively.
    Mode 3, which is not shown, is the complex conjugate of mode 2, and is also computed accurately.
    There is a larger amount of variation in the 4th and 7th modes between the five data sets, but as shown above, the shape of these modes remain consistent.
  }
  \label{fig:fhne-modes}
\end{figure}

To compute the Koopman eigenvalues, eigenfunctions, and modes, we use a polynomial kernel with $\alpha = 20$.
Furthermore, we scale each data set so that the mean value of $\|\vec x_m\| = 1$, which is equivalent to introducing a scaling parameter in the kernel $f$.
For these data sets, both $\mat{X}$ (for DMD) and $\mat{\Psi}_x$ (for the kernel method) have small singular values, so when we compute the pseudoinverse of $\mat{X}$ (for DMD) or $\mat{\Sigma}$ (for the kernel method), we retain the 150 largest singular values and set the others to zero. 
This is equivalent to projecting $\mat{X}$ and $\mat{Y}$ onto the 150 leading POD modes of $\mat{X}$ or  $\mat{\Psi}_x$ and $\mat{\Psi}_y$ onto the leading 150 kernel principal components of $\mat{\Psi}_x$ before carrying out the computation. 
In our experience, this often helps to avoid the spurious unstable eigenvalues that both DMD and Extended DMD can produce.

Figure~\ref{fig:fhne-eigenvalues} shows the approximation of the
continuous-time eigenvalues obtained for the five different data sets using this methodology.
Because we used truncated SVDs for computing pseudoinverses, both DMD and the kernel method produce eigenvalues that are, up to numerical errors, contained in the left half plane.
Furthermore, both DMD and the kernel method consistently identify the eigenvalues with $\lambda_1 \approx 0$ and $\lambda_{2,3} \approx -0.006 \pm 0.053i$ though, as shown in the figure, the variance in this pair is smaller when the kernel method is used.
Although DMD identifies other eigenvalues with more negative real parts, no other eigenvalues lie in the subset of the complex plane shown in the figure. 
The kernel method, which uses a ``richer'' set of basis functions to approximate the Koopman eigenfunctions, also identifies two additional layers of the ``pyramid'' of eigenvalues that the Koopman operator should possess~\cite{Gaspard2001}.
As a result, although the numerically computed eigenvalues are truly data dependent, the kernel approach: (i) appears to approximate the leading eigenvalues of the Koopman operator with less variance than the DMD benchmark, and (ii) identifies additional Koopman eigenvalues the DMD benchmark appears to neglect.

The $v$-part of the modes corresponding to four of the leading eigenvalues are shown in Fig.~\ref{fig:fhne-modes}.
The first mode, which has $\lambda_1\approx 0$, is a scaled version of the solution  at the fixed point where the appropriate scaling factor is determined by the corresponding Koopman eigenfunction.
The error, i.e., $\|\vec\xi_1 - \vec\xi_{\text{true}}\|$ with the normalization $\|\vec\xi\|=1$, was less than 0.01 for the five data sets used in the computation (the mean value of the error was $6.9\times 10^{-3}$).
The next pair of modes (modes 2 and 3), which have $\lambda_{2,3} \approx -0.006 \pm 0.053i$, can be shown to be the ``slowest'' pair of eigenvectors of the system linearization, and provide a useful description of the evolution near the fixed point.
The error in these modes was less than 0.02 for the five data sets (with a mean value of $1.9\times 10^{-2}$). 
As a result, not only does the kernel method consistently identify the leading Koopman modes in this problem, but it is accurate as well.

The first mode associated with \textit{nonlinear effects}, which has the  eigenvalue ${\lambda_4\approx-0.013}$, is shown in the third image of Fig.~\ref{fig:fhne-modes}.
As shown in the figure, the kernel method consistently identifies this mode, though there are differences in the mode shape that are visible to the eye if one  ``zooms in'' near $x=10$ in the figure.
In general, both DMD and the kernel method become less accurate the further ``down'' (i.e., into the left half-plane) one goes, as is demonstrated by the slight yet visible differences in the 7th mode with $\lambda_7\approx -0.019 + 0.053i$ shown in the bottom panel of Fig.~\ref{fig:fhne-modes}.
Intuitively, this is because the corresponding eigenfunctions tend to become more complicated functions; in this example, the eigenfunctions corresponding to each subsequent layer of the pyramid are higher powers of the leading two eigenfunctions.
As a result, they are less likely to lie in or near the subspace spanned by our implicitly chosen basis set, and therefore less accurately computed.

In this section, we applied the kernel approach to the FitzHugh-Nagumo PDE, and found that it consistently identified a subset of the Koopman modes and eigenvalues including those that are not associated with the system linearization.
The first three modes, where the true solutions happen to be known, are identified accurately, with a relative error of less than 2\%.
Given a sufficiently large amount of data and the proper choice of kernel, the this approach can
also identify additional ``layers'' in the pyramid of Koopman eigenvalues,
but here the variance in the numerically obtained eigenvalues and
modes is larger, which is likely due to the fact the eigenfunctions of interest to not lie entirely within the span of the basis functions used to represent them.
As a result, at a certain point the method transitions from yielding quantitatively accurate modes and eigenvalues to a more qualitative description of the underlying dynamics.

\section{Example: Experimental Data of Flow Over a Cylinder}
\label{sec:exampl-exper-data}

In this example, we apply the kernel approach to experimentally
obtained PIV data for flow over a cylinder taken at Reynolds number
413, used in previous work by \citet{Tu2014} and~\citet{Hemati2014}.
The experiment was conducted in a water channel designed to minimize
three-dimensional effects, and the computed velocity field sampled at
20 Hz at a resolution of  $135\times 80$ pixels~\cite{Tu2014}.
From these snapshots of the velocity field, we compute the vorticity via standard finite difference methods.
In the original manuscripts, a total of 8000 snapshots of the velocity field were obtained, but here we use 3000 of them in our computation.

This example is intended to highlight how we foresee the kernel approach being used in practice.
Due to the experimental nature of the data {\em and the presence of   measurement and process noise}, the  true Koopman eigenvalues and modes are not known.
However as we will show shortly, the numerically computed Koopman modes and
eigenvalues have more of the properties that the true eigenvalues should
have, which implies the kernel method is producing a more accurate
approximation than DMD.
As a result, our objective here is to demonstrate the practical
advantages of the kernel method, rather than the accuracy of the
resulting modes and eigenvalues.

\begin{figure}
	\centering
	\includegraphics{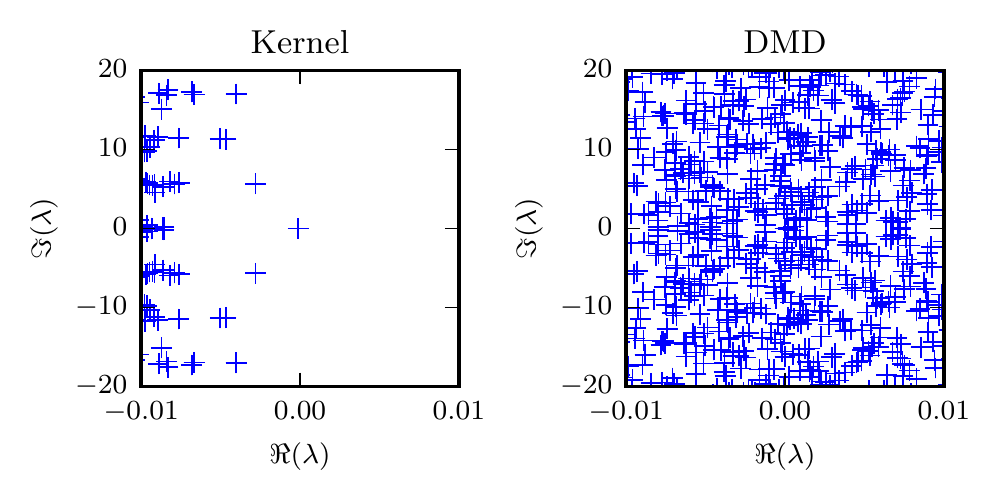}
	\caption{(left) A subset of the eigenvalues computed by applying to kernel method to a time-series of vorticity data computed from an experimentally obtained velocity field.
		(right) A subset of the eigenvalues that result from applying DMD to the same data set.
		Due to measurement and process noise in the data, the Koopman eigenvalues of this system are not known.
		However, the practical advantage of the kernel approach is that the {\em leading} eigenvalues are visually separated from the ``cloud'' of more rapidly decaying eigenvalues; while DMD also computes the eigenvalues associated with frequencies believed to exist in this system, it is not straightforward to separate these ``useful'' eigenvalues/modes from the remainder.
	}
	\label{fig:cylinder-eigenvalues}
\end{figure}

As before, we use a polynomial kernel with $\alpha = 20$, scale the
data such that the mean value of $\|\vec x_i\| = 1$, and use
(\ref{eq:kernel-edmd}) to approximate the Koopman eigenvalues and modes.
From previous efforts including \citet{Tu2014} and \citet{Hemati2014}, an effective set of modes and eigenvalues can be obtained from this
data using DMD and similar methods. 
Figure~\ref{fig:cylinder-eigenvalues} shows the leading eigenvalues obtained using the kernel method and, for reference, those computed using DMD.
In both cases, all 2999 singular values were above machine epsilon, and were retained when computing pseudoinverses.
The practical advantage of the kernel approach is that the eigenvalues and modes of interest can be identified visually; in particular, they are the ``most slowly decaying'' eigenvalues, which could be computed
efficiently even in large systems using Krylov methods.
On the other hand, DMD has a ``cloud'' of both stable and unstable eigenvalues, and the modes known to be most important  are not necessarily those
that decay the slowest.
We should note that a ``useful'' set of eigenvalues are contained in this cloud, but energy~\cite{Tu2014} or
sparsity-promoting~\cite{Jovanovic2013} methods must be used in order to select this subset of the eigenvalues.
To do this, however, both of these methods require all of eigenvalues and modes to be computed, and then truncate this ``complete'' set, which can become
computationally expensive in systems with larger numbers of snapshots and, therefore, modes to choose from.

The other eigenvalues computed by the kernel method decay more rapidly than those pictured in Fig.~\ref{fig:cylinder-eigenvalues}, and no
eigenvalues have a positive real part (up to roundoff errors).
Conceptually, the underlying system possesses a limit cycle with the addition of both process and measurement noise.
In the absence of noise, the Koopman eigenvalues would lie on the imaginary axis and in regular intervals in the left
half-plane~\cite{Gaspard2001}.
The addition of noise causes many of the eigenvalues to decay (or
decay more rapidly)~\cite{Mezic2005,Bagheri2014,WilliamsEDMD}, as the
distribution of trajectories approaches the invariant distribution that
is concentrated  near the limit cycle.
However, this decay rate depends upon the amount and type of noise added to the system, so while some decay is expected, it is unclear whether or not these eigenvalues are accurate.

\begin{figure}
	\centering
        \includegraphics{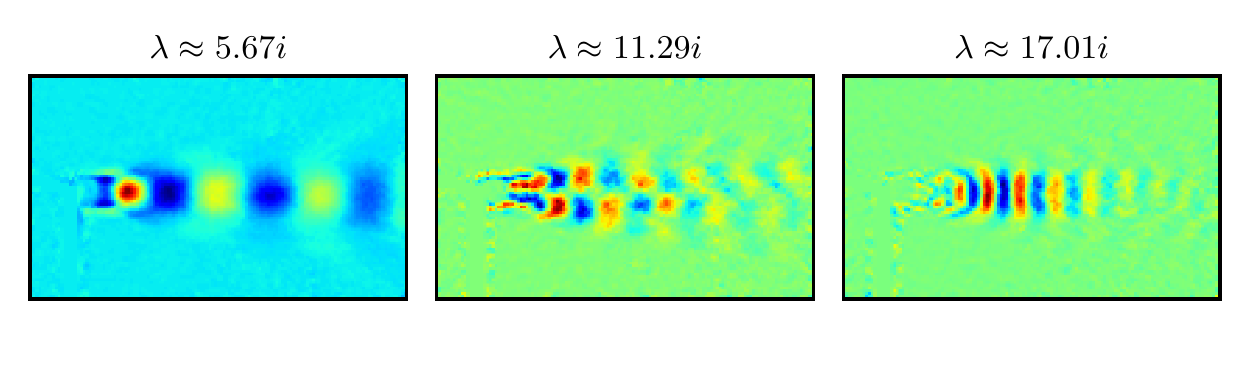}
	\caption{
          The real part of three of the leading Koopman modes and their associated eigenvalues obtained using the kernel method.
          These modes agree favorably with pre-existing modes computed using DMD~\cite{Tu2014}, but these modes {\em were identified visually} rather than via energy-based or sparsity-promoting~\cite{Jovanovic2013} methods for mode selection.
	}
	\label{fig:cylinder-modes}
\end{figure}

The real part of three of the leading Koopman modes are shown in
Fig.~\ref{fig:cylinder-modes}.
Both the shape and angular frequency (i.e., $\Im(\lambda)$) of these modes are
similar to those obtained using the full set of 8000 snapshots in
\citet{Tu2014} or \citet{Hemati2014}.
In particular, $\lambda = 5.67i$ corresponds to an oscillation
frequency of 0.90 Hz, $\lambda = 11.29i$ to a frequency of 1.79 Hz,
and $\lambda = 17.01i$ to 2.71 Hz.
In \citet{Hemati2014}, these frequencies were reported as 0.89 Hz,
1.77 Hz, and 2.73 Hz respectively. 
The modes shown in Fig.~\ref{fig:cylinder-modes} were obtained from
the left eigenvectors of the matrix $\mat{\hat{K}}$ using
(\ref{eq:mode-def}), which in this example were obtained by completely
decomposing $\mat{\hat{K}}$, but could also have been obtained by
computing the leading left-eigenvectors using the implicitly restarted
Arnoldi method~\cite{Lehoucq1998}, which would be more efficient if a
larger number of snapshots were available.

In this section, we applied the kernel method to experimental data and
computed the leading Koopman modes.
Unlike the previous example, where the kernel approach produced modes
that DMD could not, in this example both methods produce similar sets of modes.
The difference, however, is that the ``important'' modes identified by
\citet{Tu2014} and \citet{Hemati2014} using an energy-based method, are
the clearly-visible leading (i.e., most slowly-decaying) modes of the kernel method.
As a result, rather than computing {\em all\/} of the Koopman modes and
then truncate, one could compute only the {\em leading\/} Koopman modes
via (\ref{eq:mode-def}).
Furthermore, the eigenvalues identified by the kernel method possess
more of the properties that the Koopman eigenvalues should have; in
particular, they lie on the imaginary axis or the left half-plane.
Therefore, even though we cannot definitively say  that the kernel
method has accurately identified the Koopman eigenvalues, there are
conceptual and practical advantages to using the kernel method in lieu
of standard DMD even in experimental applications like the one shown here.

\section{Conclusions}
\label{sec:conclusions}

We have  presented a data-driven method for approximating the Koopman operator in
problems  with large state spaces, which commonly occurs in physical
applications such as fluid dynamics.
The kernel method we have developed defines a subspace of scalar
observables implicitly  through a kernel function, which allows the
method to approximate  the Koopman operator with the same asymptotic cost as DMD, but
with a far larger set of observables.
In this manuscript, we used a polynomial kernel, whose associated set
of basis functions can represent $\alpha$-th order polynomials.
However, many other kernels are available, associated with other basis
sets; as a result, the choice of kernel, like the
choice of a basis set in Extended DMD, is important but 
user determined.

To highlight the effectiveness of the kernel method, we applied it to
a pair of illustrative examples.
The first is the FitzHugh-Nagumo PDE in one spatial dimension, where a subset of the
Koopman eigenvalues and modes could be deduced by linearizing the system about an equilibrium point.
For this subset, the kernel method consistently and accurately identified the leading Koopman eigenvalues and modes.
However, it also accurately computed additional Koopman eigenvalues that are not associated with the system linearization and, although we do not have an analytical expression for the corresponding modes, consistently identified a mode shape.

In practical applications, the data often come from systems whose true
Koopman eigenvalues (or related information such as the location and
linearization about fixed points) are unknown.
Our second example, which used experimental data from flow over a
cylinder, was a more realistic example of how we envision the kernel
approach being applied.
In that example, the numerically obtained eigenvalues had more of the
properties the Koopman eigenvalues {\em should} have, but it was
unclear whether or not they were quantitatively accurate because the
distribution and nature of the noise in this example was unknown.
The primary advantage of the kernel method here was that the
Koopman modes of interest were associated with the eigenvalues that decay the slowest.
As a result, they could be computed using the leading left-eigenvectors
of the finite dimensional approximation, rather than requiring the
complete set of eigenvectors and eigenvalue be computed and then truncated using
energy-based or sparsity promoting methods. 

In the end, the Koopman operator is an appealing mathematical framework for
defining a set of spatial modes because these modes are intrinsic to the
dynamics rather than associated with a set of data.
However, obtaining effective approximations of this operator is
non-trivial, particularly when all that is available are data.
One method for obtaining such an approximation is Extended Dynamic
Mode Decomposition~\cite{WilliamsEDMD}, which is, for certain choices of basis functions, only computationally
feasible in problems with small state spaces.
The kernel method presented here is conceptually identical to Extended Dynamic Mode Decomposition, but computationally feasible even in large
systems; in particular, the asymptotic cost of the kernel method is identical to that of ``standard'' Dynamic Mode Decomposition, and therefore, it can be used anywhere DMD is currently being used.
Like most existing data-driven methods, there is no guarantee that the
kernel approach will produce  accurate approximations of even
the leading eigenvalues and modes, but it often appears to produce {\em
  useful} sets of modes in practice if the kernel and truncation level are chosen properly. 
As a result, approaches like the kernel method presented here are a first step
towards the ultimate goal of enabling the conceptual tools
provided by Koopman spectral analysis to be used in practical applications.

\section*{Acknowledgments} 
The authors acknowledge support from NSF (awards DMS-1204783 and
CDSE-1310173) and AFOSR.

\bibliographystyle{IEEEtranN}

\bibliography{koopman}

\end{document}